\documentclass{amc}
 \usepackage{amsmath}
  \usepackage{paralist}
  \usepackage{graphics} %% add this and next lines if pictures should be in esp format
  \usepackage{epsfig} %For pictures: screened artwork should be set up with an 85 or 100 line screen
  \usepackage{amsfonts,amssymb}
  \usepackage{float}
  \usepackage{graphicx}
  \usepackage[thinlines,thiklines]{easybmat}
  \usepackage{multirow}
  \usepackage{algorithmic}
  \usepackage{algorithm}

  \usepackage{hyperref} %% Warning: when you first run your tex file, some errors might occur, please just
\def\currentvolume{X}
 \def\currentissue{X}
  \def\currentyear{200X}
   
    \def\ppages{X--XX}

\newcommand{\Z}{\mathbb{Z}} % les entiers relatifs
\newcommand{\Q}{\mathbb{Q}} % les rationnels
\newcommand{\F}{\mathbb{F}} % un corps fini
\newcommand{\K}{\mathbb{K}} % un corps
\newcommand{\p}{\mathfrak{p}} % un p gothik
\newcommand{\Cl}{\operatorname{Cl}}

\newtheorem{theorem}{Theorem}[section]

\newtheorem{lemma}[theorem]{Lemma}

\theoremstyle{definition}
\newtheorem{definition}[theorem]{Definition}

\renewcommand{\algorithmicrequire}{\textbf{Input:}}
\renewcommand{\algorithmicensure}{\textbf{Output:}}

\title[Computation of ideal class groups]
      {Improvements in the computation of ideal class groups of imaginary quadratic number fields}

\author[Jean-Fran\c{c}ois Biasse]{}


\thanks{The author is supported by a DGA grant}

\begin{document}
\maketitle

%% Enter the first author's name and address:
\centerline{\scshape Jean-Fran\c{c}ois Biasse }
\medskip
{\footnotesize
 %% please put the address of the first author
 \centerline{LIX - \'{E}cole Polytechnique}
   \centerline{91128 Palaiseau , France}
} %% Do not forget to end the {\footnotesize by the sign }

\bigskip

%% The name of the associate editor will be entered by an editorial staff
 \centerline{(Communicated by Tanja Lange)}

\begin{abstract}
We investigate improvements to the algorithm for the computation of ideal class groups described by Jacobson in the imaginary quadratic case. These improvements rely on the large prime strategy and a new method for performing the linear algebra phase. We achieve a significant speed-up and are able to compute ideal class groups with discriminants of 110 decimal digits in less than a week. 
\end{abstract}

\section{Introduction}

%Description des particularit?s math?matiques du groupe de classes d'id?aux, notamments la d?composition suivant les id?aux premiers.
Given a fundamental discriminant $\Delta$, it is known that the corresponding ideal class group $\Cl(\Delta)$ of the order $\mathcal{O}_{\Delta}$ of discriminant $\Delta$ in $\K = \Q(\sqrt{\Delta})$ is a finite abelian group that can be decomposed as 
$$\Cl(\Delta) \simeq \bigoplus_i \Z/d_i\Z,$$
where the divisibility condition $d_{i}|d_{i+1}$ holds. 
%In cryptography we make an extensive use of finite groups in which the discrete logarithm problem (DLP) is hard to solve. Finding the structure of such groups is closely related to solving instances of the DLP and thus of both cryptographic and number theoretic interest. Some key-exchange protocols relying on the difficulty of solving the discrete logarithm problem in imaginary quadratic orders have been proposed \cite{buchmannProtocol,JacobsonProtocol}. In this paper we investigate improvements in the computation of the group structure of $\Cl(\Delta)$.\\
In this paper we investigate improvements in the computation of the group structure of $\Cl(\Delta)$: that is, determining the $d_i$, which is of both cryptographic and number theoretic interest. Indeed some cryptographic protocols relying on the difficulty of solving the discrete logarithm problem (DLP) in imaginary quadratic orders have been proposed \cite{buchmannProtocol,JacobsonProtocol}, and solving instances of the DLP is closely related to finding the group structure of $\Cl(\Delta)$.

In 1968 Shanks \cite{Shanks} proposed an algorithm relying on the baby-step giant-step method in order to compute the structure of the ideal class group of an imaginary quadratic number field in time $O\left( |\Delta |^{1/4 + \epsilon} \right)$, or $O\left( |\Delta |^{1/5 + \epsilon} \right)$ under the extended Riemann hypothesis \cite{LenstraShanks}. This allows us to compute class groups of discriminants having up to 20 or 25 decimal digits. Then a subexponential strategy was described in 1989 by Hafner and McCurley \cite{hafner}. The expected running time of this method is
%$$L\left( \sqrt{2}\right) = \left( \exp \sqrt{\log\Delta\log\log\Delta} \right)^{\sqrt{2}+o(1)}$$
$$e^{ \left( \sqrt{2}+o(1)\right) \sqrt{\log|\Delta|\log\log|\Delta|} }.$$
Buchman and D\"{u}llmann \cite{dullmann} computed class groups with discriminants of around 50 decimal digits using an implementation of this algorithm. An improvement of this method was published by Jacobson in 1999 \cite{JacobsonPhd}. He achieved a significant speed-up by using sieving strategies to generate the matrix of relations. He was able to compute the structure of class groups of discriminants having up to 90 decimal digits. More recently Sutherland \cite{Sutherland} used generic methods in order to compute class groups with discriminants having 100 decimal digits. Unlike the previous algorithms, this one relies heavily on the particular structure of $\Cl(\Delta)$ thus obtaining variable performances depending on the values of $\Delta$.

Our approach is based on that of Jacobson, using new techniques to accelerate both the sieving phase and the linear algebra phase; we have obtained the group structure of class groups of 110 decimal digit discriminants.

\section{The ideal class group}

In this section we give essential results concerning the ideal class group and the subexponential strategies for computing its structure. For a more detailed description of the theory of ideal class groups we refer to \cite{cohen} and \cite{neukirch}. In the following, $\Delta$ is a non-square integer congruent to 0 or 1 modulo 4, and the quadratic order of discriminant $\Delta$ is defined as the $\Z$-module
$$\mathcal{O}_{\Delta} = \Z + \frac{\Delta + \sqrt{\Delta}}{2}\Z.$$
We also denote by $\K$ the field $\Q(\sqrt{\Delta})$.
%In the following, $\Delta$ denotes a fundamental discriminant and $\K = \Q(\sqrt{\Delta})$ is the corresponding number field.

\subsection{Description}

%Description des particularités mathématiques du groupe de classes d'idéaux, notamments la décomposition suivant les idéaux premiers.
Elements of $\Cl(\Delta)$ are obtained from fractional ideals of $\mathcal{O}_{\Delta}$, which are $\Z$-modules of $\K$ of the form:
$$\mathfrak{a} = q \left( a\Z + \frac{b + \sqrt{\Delta}}{2}\Z\right), $$
where $a$ and $b$ are integers with $b\equiv \Delta \ \text{mod}\ 2$ and $q$ is a rational number. The prime ideals are the fractional ideals for which there exists a prime number $p$ such that:
$$\p = p\Z + \frac{b_p+\sqrt{\Delta}}{2}\Z\ \ \text{or}\ \ \p = p\Z\ (p \text{ inert in } \K).$$

\begin{definition}[Ideal Class group]
Let $\mathcal{I}_{\Delta}$ be the set of invertible fractional ideals of $\mathcal{O}_{\Delta}$, and $\mathcal{P}_{\Delta}=\left\lbrace (\alpha)\in\mathcal{I}_{\Delta},\alpha\in\K\right\rbrace $ the subset of principal ideals. We define the ideal class group of $\Delta$ as : 
$$\Cl(\Delta) := \mathcal{I}_{\Delta}/\mathcal{P}_{\Delta},$$
where the group law is the one derived from the multiplication of $\Z$-modules.
\end{definition}

For every $\mathfrak{a}\in\mathcal{I}_{\Delta}$, there exist uniquely determined prime ideals $\p_1,\hdots ,\p_n$ and exponents $e_1,\hdots , e_n$ in $\Z$ such that 
$$\mathfrak{a} = \p_1^{e_1}\hdots\p_n^{e_n}.$$
Unlike $\mathcal{I}_{\Delta}$, the ideal class group $\Cl(\Delta)$ is a finite group. Its order is called the class number and usually denoted by $h(\Delta)$. It grows like $|\Delta|^{1/2+\epsilon}$, as shown in \cite{siegel}.

\subsection{Computing the group structure}
%Description haut niveau de la stratégie de recherche de relations ainsi que des phases d'algèbre linéaires qui s'ensuivent.
The algorithm for computing the group structure of $\Cl(\Delta)$ is divided into two major phases: relation collection and linear algebra. In the first phase, we begin by precomputing a factor base $\mathcal{B} = \left\lbrace \p_1,\hdots,\p_n\right\rbrace $ of non-inert prime ideals satisfying $\mathcal{N}\left( \p_i\right) \leq B$, where $B$ is a smoothness bound. Then we look for relations of the form
$$\left( \alpha\right) = \p_1^{e_1}\hdots\p_n^{e_n}, $$
where $\alpha\in\K$. Every $n$-tuple $[e_1,\hdots,e_n]$ collected becomes a row of what we will refer to as the relation matrix $A\in\Z^{m\times n}$. We have from \cite{bach} the following important result:
\begin{theorem}
Let $\Lambda$ be the lattice spanned by the set of the possible relations. Assuming GRH, if $B\geq 6\log^2\Delta$, then we have
$$\Cl(\Delta) \simeq \Z^n/\Lambda.$$
\end{theorem}

After the relation collection phase we can test if $A$ has full rank and if its rows generate $\Lambda$ using methods described in \textsection \ref{hnf_algo}. If it is not the case then we have to compute more relations. From now on we assume that $A$ has full rank and that its rows generate $\Lambda$.

The linear algebra phase consists of computing the Smith Normal Form (SNF) of $A$. Any matrix $A$ in $\Z^{n\times n}$ with non zero determinant can be written as
\[ A = V^{-1}\left( \begin{array}{cccc}
d_1    & 0      & \hdots & 0      \\
0      & d_2    & \ddots & \vdots \\
\vdots  & \ddots & \ddots & 0 \\
0      & \hdots & 0      & d_n   \end{array} \right) U^{-1}\], 
where $d_{i+1} | d_i$ for all $1 \leq i < n$ and $U$ and $V$ are unimodular matrices in $\Z^{n\times n}$. The matrix $\text{diag}(d_1,\hdots,d_n)$ is called the SNF of $A$. If $m=n$ and $\text{diag}(d_1,\hdots,d_n) = \text{SNF}(A)$ then 
$$\Cl(\Delta) \simeq \bigoplus_{i=1}^{n} \Z/d_i\Z.$$
This reduces the problem of computing the group structure of $\Cl(\Delta)$ to computing the SNF of a relation matrix $A$ in $\Z^{n\times n}$. For an arbitrary $A$ in $\Z^{m\times n}$ we start by computing the Hermite Normal Form (HNF) of $A$. A matrix $H$ is said to be in HNF if it has the shape
%\[ H = \left( \begin{array}{cccc}
%h_{1,1}& 0      & \hdots & 0      \\
%\vdots & h_{2,2}& \ddots & \vdots \\
%\vdots & \vdots & \ddots & 0      \\
%*      & *      & \hdots & h_{n,n}\\
%   0   & 0      & \hdots      & 0 \\
%\end{array} \right)\] 
%\[ H = \left( 
%   \begin{BMAT}(e)[2pt,3cm,3cm]{cccc}{cccc.c}
%h_{1,1}& 0      & \hdots & 0      \\
%\vdots & h_{2,2}& \ddots & \vdots \\
%\vdots & \vdots & \ddots & 0      \\
%*      & *      & \hdots & h_{n,n}\\
%   0   & \hdots      & \hdots      & 0 
%%  \begin{BMAT}{c}{c} & & & 
%%	(0)
% % \end{BMAT}
%\end{BMAT}
%   \right)  \]
\[ H = \left( 
   \begin{BMAT}(@)[2pt,3cm,3cm]{c}{c.c}
   \begin{BMAT}(e){cccc}{cccc}
h_{1,1}& 0      & \hdots & 0      \\
\vdots & h_{2,2}& \ddots & \vdots \\
\vdots & \vdots & \ddots & 0      \\
*      & *      & \hdots & h_{n,n} 
  \end{BMAT} \\
\begin{BMAT}[2pt,3cm,1cm]{c}{c} 
	(0)
\end{BMAT}
\end{BMAT}
   \right)  \],
where $0\leq h_{ij} < h_{ii}$ for all $j<i$ and $h_{ij}=0$ for all $j>i$. For each matrix $A$ in $\Z^{m\times n}$ there exists a matrix $H$ in HNF and a unimodular matrix $W$ in $\Z^{m\times m}$ such that
$$H = WA.$$
The upper block of $H$ is a $n\times n$ relation matrix whose SNF provides us the group structure of $\Cl(\Delta)$. There is an index $l$ such that $h_{i,i} = 1$ for every $i\geq l$. The upper left $l\times l$ submatrix of $H$ is called the essential part of $H$. In order to compute the group structure of $\Cl(\Delta)$ it suffices to compute the SNF of the essential part of $H$, which happens to have small dimension in our context.
\subsection{The use of sieving for computing the relation matrix}
%Description de la creation de relations via sieving. J'attends un chapitre du nouveau livre de Jacobson sur l'équation de Pell qui d'après lui apporte un éclairage nouveau sur la relation idéaux-formes quadratiques.
The use of sieving to create the relation matrix was first described by Jacobson \cite{JacobsonPhd}. Here we follow the approach of \cite{JacobsonPell} Chap.13, which relies on the following lemma:

\begin{lemma}
 If $\mathfrak{a} =  \left( a\Z + \frac{b + \sqrt{\Delta}}{2}\Z\right)$ with $a>0$, then for all $x,y$ in $\Z$ there exists $\mathfrak{b}\in\mathcal{I}_{\Delta}$ such that $\mathfrak{a}\mathfrak{b}\in\mathcal{P}$ and
$$\mathcal{N}(\mathfrak{b}) = ax^2 + bxy + \frac{b^2 - \Delta}{4a}y^2.$$
\end{lemma}
%\begin{lemma}
%If $\gamma\in\mathfrak{a} = a\Z + \frac{b + \sqrt{\Delta}}{2}\Z$ with $a>0$ and $x,y\in\Z$ are such that
%$\gamma =  ax + \frac{b+\sqrt{\Delta}}{2}y$
%then there exists an ideal $\mathfrak{b}\in\mathcal{I}$ satisfying $(\gamma)=\mathfrak{a}\mathfrak{b}$ and 
%$$\mathcal{N}(\mathfrak{b}) = ax^2 + bxy + \frac{b^2 - \Delta}{4a}y^2$$
%\end{lemma}

%Here we restricted ourselves to fractional ideals $\mathfrak{a}$ of the form 
%$$\mathfrak{a} = a\Z + (b + \sqrt{\Delta})/2\Z$$
%with $a>0$. We do not loose generality because if $\mathfrak{a} = q \left( a\Z + (b + \sqrt{\Delta})/2\Z\right)$ then its ideal class in %$Cl(\Delta)$ is the same as the ideal class of 
%$$\mathfrak{a'} = |a|\Z + (b + \sqrt{\Delta})/2\Z$$
The strategy for finding relations is the following: We start with 
$$\mathfrak{a}=\prod_i \p_i^{e_i} =: \left( a\Z + \frac{b + \sqrt{\Delta}}{2}\Z\right),$$
whose norm is $B$-smooth. Then we choose a sieve radius $R$ satisfying $R\approx \sqrt{|\Delta|/2}/\mathcal{N}(\mathfrak{a})$ and we look for values of $x\in[-R,R]$ such that $\varphi(x,1)$ is $B$-smooth where
$$\varphi(x,y) = ax^2 + bxy + \frac{b^2 - \Delta}{4a}y^2,$$
which allows us to find $\mathfrak{b} = \prod_i\p_i^{f_i}$ satisfying $\mathfrak{a}\mathfrak{b}=(\gamma)$ for some $\gamma$ in $\K$. The $\p_i$ and $f_i$ are deduced from the decomposition $\varphi(x,1)=\prod_ip_i^{v_i}$. For more details we refer to \cite{JacobsonPell}, Chap 13. This method yields the relation 
$$(\gamma) = \prod_i \p_i^{e_i+f_i}.$$
Now given a binary quadratic form $\varphi(x,y)=ax^2+bxy+cy^2$ of discriminant $\Delta$, we are interested in finding values of $x\in[-R,R]$ such that $\varphi(x,1)$ is $B$-smooth. This can be done trivially by testing all the possible values of $x$, but there is a well-known method for pre-selecting some values of $x$ in $[-R,R]$ that are going to be tested, namely the quadratic sieve (introduced by Pomerance \cite{pomerance82}). It consists in initializing to 0 an array $S$ of length $2R+1$ and precomputing the roots $r_i'$ and $r_i''$, or the double root $r_i'$,  of $\varphi(x,1)\mod p_i$ for each $p_i\leq B$ such that $\left( \frac{\Delta}{p_i} \right) \neq -1$ . Then for each $x$ in $[-R,R]$ of the form $x=r_i+kp_i$ for some $k$, we add $\lfloor\log p_i\rfloor$ to $S[x]$. At the end of this procedure, if $\varphi(x,1)$ is $B$-smooth, then $S[x]\approx\log\varphi(x,1)$. As $\varphi(x,1)\approx \sqrt{\Delta/2}R$, we set a bound 
\begin{equation}\label{eqK}
F = \log\left( \sqrt{\frac{\Delta}{2}}R\right) -T\log(p_n), 
\end{equation}
%$$K = \log\left( \sqrt{\frac{\Delta}{2}}R\right) -T\log(p_n)$$ 
where $T$ is a number representing the tolerance to rounding errors due to integer approximations. We then perform a trial division test on every $\varphi(x,1)$ such that $S[x]\geq F$.

\section{Practical improvements}

In this section we describe the improvements that allowed us to achieve a significant speed-up with respect to the existing algorithm and the computation of class group structures of large discriminants. Our contribution is to take advantage of the large prime variants, of an algorithm due to Vollmer \cite{vollmer} for the SNF which had not been implemented in the past, and of special Gaussian elimination techniques.

\subsection{Large prime variants}

%Quid des large primes dans différents contextes. Description théorique pour un , deux ou plus de large primes.
The large prime variants were developed in the context of integer factorization to speed up the relation collection phase in both the quadratic sieve and the number field sieve. Jacobson considered analogous variants for class group computation \cite{JacobsonPhd}, but the speed-up of the relation collection phase was achieved at the price of such a slow-down of the linear algebra that it did not significantly improve the overall time. The main idea is the following: We define the ``small primes" to be the prime ideals in the factor base and the small prime bound as the corresponding bound $B_1=B$. Then we define a large prime bound $B_2$. During the relation collection phase we choose not to restrict ourselves to relations only involving primes $\p$ in $\mathcal{B}$ but we also keep relations of the form 
$$(\alpha)=\p_1\hdots\p_n \p \ \ \text{and}\ \ (\alpha)=\p_1\hdots\p_n \p\p'$$
for $\p_i$ in $\mathcal{B}$, and for $\p,\p'$ of norm less than $B_2$. We will respectively refer to them as 1-partial relations and 2-partial relations. Keeping partial relations only involving one large prime is the single large prime variant, whereas keeping two of them is the double large prime variant which was first described by Lenstra and Manasse \cite{Lenstra2LP}. In this paper we do not consider the case of more large primes, but it is a possibility that has been studied in the context of factorization \cite{Lenstra3LP}.

Partial relations may be identified as follows. Let $m$ be the residue of $\varphi(x,1)$ after the division by all primes $p\leq B_1$, and assume that $B_2 < B_1^2$. If $m=1$ then we have a full relation. If $m\leq B_2$ then we have a 1-partial relation. We can see here that detecting 1-partial relations is almost for free. If we also intend to collect 2-partial relations then we have to consider the following possibilities:
\begin{enumerate}
 \item $m > B_2^2$;
 \item $m$ is prime and $m > B_2$;
 \item $m \leq B_2$;
 \item $m$ is composite and $B_1^2 < m \leq B_2^2$.
\end{enumerate}
In Cases 1 and 2 we discard the relation. In Case 3 we have a 1-partial relation, and in Case 4 we have $m=pp'$ where $p = \mathcal{N}(\p)$ and $p' = \mathcal{N}(\p')$. After testing if we are in Cases 1, 2, or 3 we have to factorize the residue. We have done that using Milan's implementation of the SQUFOF algorithm \cite{tifa} based on the theoretical work of \cite{squfof}.

Even though we might have to factor the residue, collecting a partial relation is much faster than collecting a full relation because the probability that $\mathcal{N}(\mathfrak{b})$ is $B_2$-smooth is much greater than the probability that it is $B_1$-smooth. This improvement in the speed of the relation collection phase comes at a price: The number of columns in the relation matrix is much greater, thus preventing us from running the linear algebra phase directly on the resulting relation matrix and forcing us to find many more relations since we have to produce a full rank matrix. We will see in \textsection \ref{gauss} how to reduce the dimensions of the relation matrix using Gaussian elimination techniques and in \textsection \ref{opt} how to optimize the parameters to make the creation of the relation matrix faster, even though there are many more relations to be found.

\subsection{Gaussian elimination techniques}\label{gauss}
%Une section qui risque de devoir comporter des schémas explicatifs. Description générale de la méthode par calcul d'un arbre couvrant de poid minimal, importance de la fonction de coup ainsi que de la suppression régulière de lignes dont les coefficients sont trop gros.
Traditionally rows were recombined to give full relations as follows: In the case of 1-partial relations, any pair of relations involving the same large prime $\p$ were recombined into a full relation. In the case of 2-partial relations, Lenstra \cite{Lenstra2LP} described the construction of a graph whose vertices were the relations and whose edges linked vertices having one large prime in common. Finding independent cycles in this graph allows us to find recombinations of partial relations into full relations.

In this paper we rather follow the approach of Cavallar \cite{Cavallar}, developed for the number field sieve, which uses Gaussian elimination on columns without distinguishing those corresponding to the large primes from the others. One of the main differences between our relation matrices and the matrices produced in the number field sieve is that our entries are in $\Z$ rather than $\F_2$, thus obliging us to monitor the evolution of the size of the coefficients. Indeed, eliminating columns at the price of an explosion of the size of the coefficients can be counter-productive in preparation for the HNF algorithm.

In what follows we will use a few standard definitions that we briefly recall here. First, subtracting two rows is called \textit{merging}. This is because rows are stored as lists of the non-zero entries sorted with respect to the corresponding columns and subtracting them corresponds to merging the two sorted lists. If two rows $r_1$ and $r_2$ share the same prime $\p$ with coefficients $c_1$ and $c_2$ respectively then multipling $r_1$ by $c_2$ and $r_2$ by $c_1$ and merging is called \textit{pivoting}. %In naive gaussian elimination strategies, when eliminating a given column, we chose a pivot at random and use it in order to eliminate the prime $\p$ that it shares with the other rows of the column without optimizing the resulting growth of the density of the matrix and of he size of the coefficients. Finally, given a prime ideal $\p$, eliminating every contribution of $\p$ in the rows contained in the corresponding column except for the last pivot using merges is called a $k$-way merge if $k$ is the weight of the column. The last pivot is then deleted thus eliminating the column.\\
Finally, finding a sequence of pivots leading to the elimination of a column of Hamming weight $k$ is a $k$-way merge.

We aim to reduce the dimension of the relation matrix by performing $k$-way merges on the columns of weight $k=1,\hdots,w$ in increasing order for a certain bound $w$. Unfortunately, the density of the rows and the size of the coefficients increase during the course of the algorithm, thus obliging us to use optimized pivoting strategies. In what follows we describe an algorithm performing $k$-way merges to minimize the growth of both the density and the size of the coefficients.

First we have to define a cost function defined over the set of the rows encapsulating the difficulty induced for the HNF algorithm. In factorization, we want to find a vector in the kernel of the relation matrix which is defined over $\F_2$; the only property of the row that really matters is its Hamming weight. In our context, we need to minimize the Hamming weight of the row, but we also have to take into account the size of the coefficients. Different cost functions lead to different elimination strategies. Our cost function was determined empirically: We took the number of non-zero entries, counting $c$ times those whose absolute value was above a bound $Q$, where $c$ is a positive number. If $r = [e_1,\hdots,e_n]$ corresponds to $(\alpha)=\prod_i\p_i^{e_i}$ then
$$C(r) = \sum_{1\leq|e_i|\leq Q}1 + c\sum_{|e_j| > Q}1.$$
Indeed as we will see, matrices with small entries are better suited for the HNF algorithm described in \textsection \ref{hnf_algo}. Let us assume now that we are to perform a $k$-way merge on a given column. We construct a complete graph $\mathcal{G}$ of size $k$ as follows:
\begin{itemize}
 \item The vertices are the rows $r_i$.
 \item Every edge linking $r_i$ and $r_j$ is labeled by $C(r_{ij})$, where $r_{ij}$ is obtained by pivoting $r_i$ and $r_j$.
\end{itemize}
Finding the best sequence of pivots with respect to the cost function $c$ we chose is equivalent to finding the minimum spanning tree $\mathcal{T}$ of $\mathcal{G}$, and then recombining every row $r$ with its parent starting with the leaves of $\mathcal{T}$.

Unfortunately, some coefficients might grow during the course of column eliminations despite the use of this strategy. Once a big coefficient is created in a given row $r$, it is likely to spread to other rows once $r$ is involved in another column elimination. We must therefore discard such rows as quickly as possible. In our implementation we chose to do it regularly: Once we have performed all the $k$-way merges for  $k\leq 10\cdot i$ and $i=1,\hdots,w/10$ we discard a fixed number $K$ of the rows containing the largest coefficients.

We show in Table \ref{TabCrunch} the effect of the use of a cost function taking into account the size of the coefficients and the regular discard of the worst rows for $\Delta = -4(10^{70}+1)$ with $c = 100$, $Q = 8$ and $K = 10$. We kept track of the evolution of the dimensions of the matrix, the average Hamming weight of the rows, and the maximum and minimum size of the coefficients. In the first case we use the traditional cost function that only takes into account the Hamming weight of the rows and we keep deleting the worst rows regularly; this corresponds to taking $c=1$ and $K=10$. In the second case, we use the cost function described above but without row elimination by setting $c=100$ and $K=0$. In the third case, we combine the two ($c=100$ and $K=10$). We clearly see that the coefficients are properly monitored only in the latter case. Indeed using a cost function that does not take into account the size of the coefficients and just discarding the worst rows regularly seems more efficient in terms of reduction of the matrix dimension, but the row corresponding to $i=12$ (that is to say after all the 120-way merges) clearly shows that we run the risk of an explosion of the coefficients.

\begin{figure}[h!]
\caption{Comparative table of elimination strategies}
\label{TabCrunch}
\small 
\begin{center}
 \begin{tabular}{|c|c|c|c|c|c|}
\hline
  \multicolumn{6}{|c|}{\textbf{Without score depending on the size of the coefficients}} \\
\hline
 $i$ & Row Nb & Col Nb & Average weight & max & min \\
\hline
0 & 38752 & 45975 & 22 & 10 & -10\\
%1 & 2572 & 1896 & 53.3165 & 14 & -14\\
2 & 2334 & 1668 & 76 & 21 & -20\\
%3 & 2203 & 1547 & 97.8225 & 31 & -31\\
4 & 2123 & 1477 & 117 & 52 & -56\\
%5 & 2065 & 1429 & 135.486 & 50 & -53\\
6 & 2028 & 1402 & 146 & 59 & -62\\
%7 & 1982 & 1366 & 163.048 & 80 & -79\\
8 & 1951 & 1345 & 175 & 72 & -65\\
%9 & 1916 & 1320 & 191.673 & 154 & -114\\
10 & 1890 & 1304 & 203 & 193 & -196\\
%11 & 1862 & 1286 & 215.907 & 364 & -408\\
12 & 1836 & 1270 & 219 & 212 & -2147483648\\
\hline
\multicolumn{6}{|c|}{\textbf{Without row elimination}} \\
\hline
 $i$ & Row Nb & Col Nb & Average weight & max & min \\
\hline
0 & 38752 & 45975 & 22 & 10 & -10\\
%1 & 2596 & 1910 & 53.6086 & 16 & -18\\
2 & 2373 & 1687 & 79 & 30 & -40\\
%3 & 2277 & 1591 & 99.7159 & 53 & -38\\
4 & 2224 & 1538 & 118 & 67 & -50\\
%5 & 2192 & 1506 & 131.286 & 88 & -87\\
6 & 2158 & 1472 & 148 & 71 & -132\\
%7 & 2133 & 1447 & 164.724 & 2584 & -9672\\
8 & 2117 & 1431 & 179 & 2648 & -10568\\
%9 & 2109 & 1423 & 185.855 & 193536 & -282624\\
10 & 2097 & 1411 & 196 & 347136 & -337920\\
%11 & 2087 & 1401 & 207.378 & 265814016 & -233963520\\
12 & 2080 & 1394 & 214 & 268763136 & -173162496\\
\hline
\multicolumn{6}{|c|}{\textbf{With adapted score and row elimination}} \\
\hline
 $i$ & Row Nb & Col Nb & Average weight & max & min \\
\hline
0 & 38752 & 45975 & 22 & 10 & -10\\
%1 & 2586 & 1910 & 53.3333 & 13 & -13\\
2 & 2357 & 1691 & 76 & 17 & -17\\
%3 & 2246 & 1590 & 96.9145 & 21 & -22\\
4 & 2176 & 1530 & 114 & 27 & -30\\
%5 & 2118 & 1482 & 131.667 & 36 & -32\\
6 & 2074 & 1448 & 149 & 37 & -37\\
%7 & 2037 & 1421 & 167.71 & 47 & -48\\
8 & 2013 & 1407 & 177 & 43 & -43\\
%9 & 1983 & 1387 & 188.774 & 47 & -46\\
10 & 1958 & 1372 & 199 & 44 & -45\\
%11 & 1927 & 1351 & 215.099 & 51 & -51\\
12 & 1908 & 1342 & 224 & 54 & -53\\
\hline
 \end{tabular}
\end{center}
\end{figure}
\normalsize

\subsection{Vollmer's algorithm for computing the HNF}\label{hnf_algo}

%Description de l'algorithme de Vollmer, et explications sur l'implantation, l'utilisation d'IML, les possibilités de linbox. ...
In \cite{JacobsonPhd} it has been observed that the algorithm used to compute the HNF of the relation matrix relied heavily on the sparsity of the matrix. While recombinations of the kind described in \cite{Lenstra2LP} or the techniques of \textsection\ref{gauss} reduce the dimensions of the matrix, they also dramatically increase the density of the matrix, thus slowing down the computation of the HNF. We had to find an HNF algorithm whose features were adapted to our situation. Vollmer described in \cite{vollmer} an algorithm of polynomial complexity depending on the capacity to solve diophantine linear systems, but not on the density of the matrix. It was not implemented at the time because there was no efficient diophantine linear system solver available. We implemented Vollmer's algorithm using the IML \cite{iml} library provided by Storjohann.

Here we give a brief description of the algorithm (for more details we refer to \cite{vollmer}). We assume we have created an $m\times n$ relation matrix $A$ of full rank. For each $i\leq n$, we define two matrices 
\[ A_i = \left( 
   \begin{BMAT}(e)[2pt,3cm,3cm]{ccc}{ccc}
      a_{1,1} & \hdots & a_{m,1}\\
	\vdots &  & \vdots \\ 
	a_{1,i} & \hdots & a_{m,i}
   \end{BMAT}
   \right)  \ \ \text{and}\ \ e_i = \left( 
   \begin{BMAT}(e)[2pt,0pt,3cm]{c}{cccc}
     0 \\ \vdots \\ 0 \\ 1
   \end{BMAT}
   \right).  \]
For each $i$, we define $h_i$ to be the minimal denominator of a rational solution of the system
$$A_ix = e_i;$$
this is computed using the function \texttt{MinCertifiedSol} of IML, which is an implementation of (Special)MinimalSolution from \cite{mulders}, and used in \cite{vollmer} for the complexity analysis. In \cite{vollmer} it is shown that 
$$h(\Delta) = \prod_i h_i.$$
Fortunately, analytic formulae allow us to compute a bound $h_*$ such that 
$$h_*\leq h(\Delta) < 2h_*,$$ 
so we do not have to compute $h_i$ for every $i\in[1,n]$. In addition, the matrices produced for the computation of the group structure of $\Cl(\Delta)$ have small essential part, which keeps the number of diophantine systems to solve small (about the same size as the number of columns of the essential part) as shown in \cite{vollmer}.
\begin{algorithm}[H]
\caption{Computation of the class number}
\begin{algorithmic}
\REQUIRE $\Delta$, relation matrix $A$ of full rank and $h_*$
\ENSURE $h (\Delta) $
\STATE $h\leftarrow 1$
\STATE $i\leftarrow n$
\WHILE {$h < h_*$}
\STATE Compute the minimal denominator $h_i$ of a solution of $A_i\cdot x=e_i$
\STATE $h\leftarrow h\cdot h_i$
\STATE $i\leftarrow i-1$ 
\ENDWHILE
\RETURN $h$
\end{algorithmic}
\end{algorithm}

We can compute the essential part of the HNF of $A$ with a little extra effort involving only modular reductions of coefficients; we refer to \cite{vollmer} for more details. This part of the algorithm is highly dependent on the performance of the diophantine solver we use, which in turn is mostly influenced by the number of columns of the matrix and the size of the coefficients. The bechmarks available \cite{iml} show that the algorithm runs much faster on matrices with 3-bit coefficients, which is why we took coefficient size into account in the cost function for the Gaussian elimination.

\section{Optimization of the parameters}\label{opt}

%Section montrant l'importance de l'optimisation de certains paramètres. Comportera pas mal de tableaux, peut etre des graphiques si cela est pertinent. 
In this section we proceed to optimize the parameters involved in the relation collection phase. Each parameter has an effect on the overall time taken to compute the group structure of $\Cl(\Delta)$. Recall  \eqref{eqK} giving the bound $F$; when we collect partial relations it should be adapted in the following way:
$$F = \log\left( \sqrt{\frac{\Delta}{2}} R \right) -T\log B_2,$$
where $B_2$ is the large prime bound.

\subsection{Optimization of $T$}

\begin{figure}[!h]
\caption{Optimum value of $T$}
\label{opt_T}
\begin{center}
\includegraphics[angle=-90,scale=0.3]{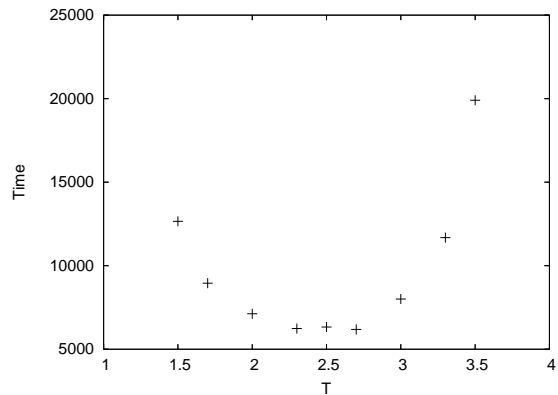}
\end{center}
\end{figure}

The parameter $T$ represents the tolerance to rounding errors in the traditional sieving algorithms. Its value is empirically determined, and usually lies in the interval $[1,2]$. In the large prime variant it also encapsulates the number of large primes we want to allow. Indeed, if there were no rounding errors one would expect this value to be 1 for one large prime and 2 for two large primes. In practice, we can exhibit an optimum value which differs slightly from what we would expect. In figure \ref{opt_T} we show the overall running time of the algorithm when the parameter $T$ varies between 1.5 and 3.5 for the discriminant $\Delta = -4(10^{75}+1)$. The size of the factor base taken is 3250, the ratio $B_2/B_1$ equals 120, and we allow two large primes.

\begin{figure}[!h]
\caption{Effect of $|\mathcal{B}|$ on the optimal value of $T$}
\label{opt_FB}
\begin{center}
\includegraphics[angle=-90,scale=0.3]{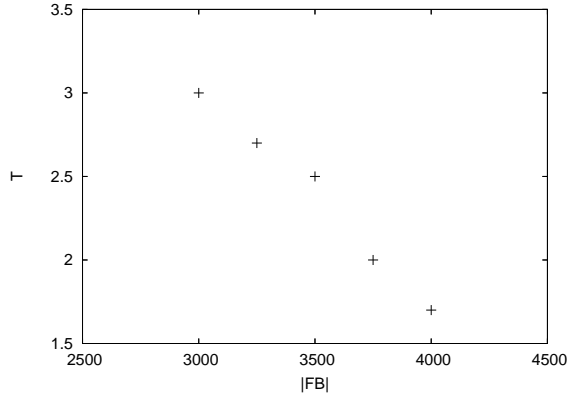}
\end{center}
\end{figure}

One of the main issues for determining the optimal value of $T$ is that it tends to shift when one modifies the value of $B_1$, the rest being unchanged. Indeed, if for example $B_2/B_1 = 120$ then 
$$F = \log\left( \sqrt{\frac{\Delta}{2}} R \right) -T\log 120B_1,$$
so when we increase $B_1$ we have to lower $T$ to compensate. Figure \ref{opt_FB} illustrates this phenomenon on the example $\Delta = -4(10^{75}+1)$, with two large primes.

\begin{figure}[!h]
\caption{Optimal value of $T$ when $n$ varies}
\label{opt_FB_n}
\begin{center}
\includegraphics[angle=-90,scale=0.3]{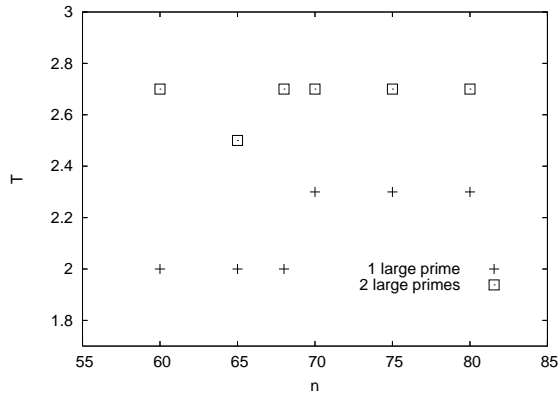}
\end{center}
\end{figure}

In Figure \ref{opt_FB_n} we study the evolution of the optimal value of $T$ for the single and double large prime variants on discriminants of the form $-4(10^n+1)$ where $n$ ranges between 60 and 80. It appears that, as we expected, the optimal value for the double large prime variant is greater than the one corresponding to the single large prime variant. This value is between 2 and 2.3 for one large prime and around 2.7 when we allow two large primes.

\subsection{ The size of the factor base}

\begin{figure}[!h]
\caption{Optimal value of $|\mathcal{B}|$}
\label{opt_FB_time}
\begin{center}
\includegraphics[angle=-90,scale=0.3]{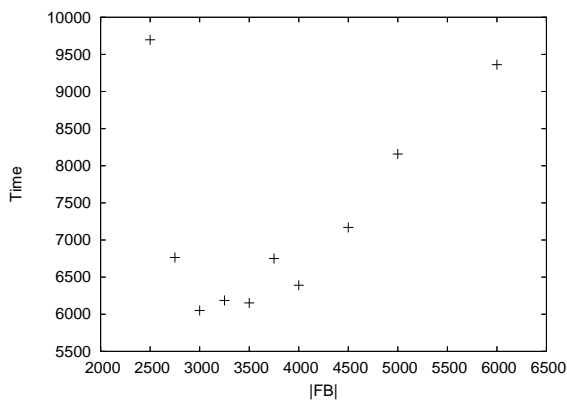}
\end{center}
\end{figure}

The optimal size of the factor base reflects the trade-off between the time spent on the relation collection phase and on the linear algebra phase. This optimum is usually not the size that minimizes the time spent on the relation collection phase. To illustrate this, Figure \ref{opt_FB_time} shows the time taken by the algorithm for $\Delta = -4(10^{75}+1)$ with $B_2/B_1 = 120$ and the corresponding optimal $T$.

\begin{figure}[!h]
\caption{Optimal value of $|\mathcal{B}|$ when $n$ varies}
\label{opt_FB_time_n}
\begin{center}
\includegraphics[angle=-90,scale=0.3]{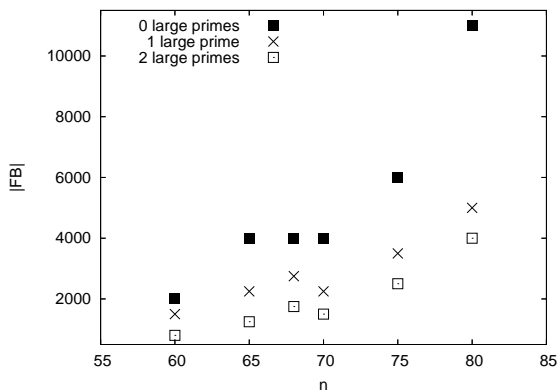}
\end{center}
\end{figure}

The optimal size of the factor base increases with the size of the discriminant. Figure \ref{opt_FB_time_n} shows the optimal size of the factor base for discriminants of the form $-4(10^n+1)$ as $n$ ranges between 60 and 80 for both one large prime and two large primes. We notice that the single large prime variant requires smaller factor bases than without large primes, and bigger factor bases than the double large prime variant. 

%We notice that for a given $\Delta = -4(10^n+1)$ the single large prime variant requires a bigger factor base than the double large prime variant. It shows that we managed to decrease the size of the factor base, thus reducing the dimension of the relation matrix after the Gaussian elimination.

\subsection { The ratio $B_2/B_1$ }

Theoretically $B_2$ should not exceed $B_1^2$. In practice, when the ratio $B_2/B_1$ is too high we lose time taking into account partial relations involving primes that are so large that they are very unlikely to occur twice and to lead to a recombination. This phenomenon is known in the context of factorization, and 120 is a common choice of value of $B_2/B_1$ (see \cite{contini}). We ran experiments using 12, 120 and 1200 as values for the ratio $B_2/B_1$. Figure \ref{TabRatio} shows the results for $\Delta = -4(10^{75}+1)$ with two large primes. We give the optimum timings for each value of the size of the factor base, and compare those values for the three different ratios. It appears that 120 is indeed the best choice, but the performance of the algorithm is not highly dependent on this parameter.

\begin{figure}[h!]
\caption{Comparative table of $B_2/B_1$}
\label{TabRatio}
%\small 
\begin{center}
 \begin{tabular}{|c|c|c|c|}
\hline
$|FB|$ & 12 & 120 & 1200 \\
\hline
3000 &  6399.60 & \textbf{6051.11} & 6173.66 \\
3250 &  6795.43 & \textbf{6185.67} & 6754.02 \\
3500 &  \textbf{6539.69} & 6821.77 & 6754.02 \\
3750 &  6916.93 & \textbf{6750.88} & 7456.92 \\
4000 &  6671.18 & \textbf{6390.48} & 7009.72\\
\hline
 \end{tabular}
\end{center}
\end{figure}

\section{Computational results}

\subsection{Comparative timings}

In Figure \ref{TabComp} we give comparative timings in seconds between no large primes and the large prime variants for discriminants of the form $-4(10^n+1)$, for $n$ between 60 and 80. We used 2.4GHz Opterons with 16GB of memory, and the NTL library with GMP. It appears that we achieved a significant speed-up by using the large prime strategy. Direct comparison with previous methods based on sieving is hard since the timings available in \cite{JacobsonPhd} were obtained on 296 MHz UtraSPARC-II processors; therefore we just quote that the computation of the group structure corresponding to $\Delta = -4(10^{80}+1)$ took 5.37 days (463968 CPU time) at the time. We also notice that the double large prime variant does not provide an impressive improvement on the overall time for the sizes of discriminant involved. The performance is comparable for discriminants of 60 decimal digits and starts showing an improvement when we reach 75 digit discrimants. %It suggests that the usual sizes of discriminant do not require the double large prime variants which is more of interest for breaking records. 

\begin{figure}[h!]
\caption{Comparative table of the performances (CPU time)}
\label{TabComp}
%\small 
\begin{center}
 \begin{tabular}{|c|c|c|c|}
\hline
$n$ & 0 Large primes & 1 Large prime & 2 Large primes \\
\hline
60 &  374 &  284& 280 \\
65 &  1019 & 756 &  776\\
68 &  2010 & 1489 & 1122 \\
70 &  2148 & 1663& 1680 \\
75 &  8409 & 6669 & 5347 \\
80 &  21215 & 17123 & 14664\\
\hline
 \end{tabular}
\end{center}
\end{figure}

\subsection{Large discriminants}

%des records pour en mettre plein la vue.
In the imaginary case, the largest class groups that had been computed using relation collection methods had 90 digits; some 100 decimal digit discriminant class group structures could be computed using the techniques of \cite{Sutherland}. With the techniques described in this paper, we achieved the computation of a class group with a 110 decimal digit discriminant. We used 100 Core2 Duo 2.4GHz Pentium IV processors with 2 GB of memory each for the sieving, and one 2.66 GHz Opteron with 64 GB of memory for the linear algebra, which is the real bottleneck of this algorithm. Indeed, the sieving phase can be trivially parallelized for as many processors as we have and does not require much memory, whereas the linear algebra can only be parallelized into the number of factors of $h$ that we get from Vollmer's algorithm (around 10 in our examples) and requires a lot of memory. Indeed the limit in terms of matrix dimensions for the diophantine solver on a 64GB memory computer seems to be around 10000 columns. For comparison, in the case of the 110 decimal digit discriminant we had to handle an 8000-column matrix (after the Gaussian reduction).

%\begin{figure}[h!]
%\caption{$h(\Delta)$ for $\Delta = -4(10^n+1)$}
%\label{Tab_h}  
%\begin{center}
%\begin{tabular}{|l|l|}
%\hline
%n & h \\
%\hline
%100 & 1462491779472195274571694315857495335176880879072 \\
%110 & 17564474658714199841017764109214000820806234521746309120\\
%\hline
%n & decomposition \\
%\hline
%100 & [2 2 2 2 2 2 2 *]\\
%110 & [2 2 2 2 2 2 2 2 2 2 2 *]\\
%\hline
%\end{tabular}
%\end{center}
%\end{figure}

\begin{figure}[h!]
\caption{Decomposition of $\Cl(\Delta)$ for $\Delta = -4(10^n+1)$}
\label{Tab_decomp}  
\begin{center}
\begin{tabular}{|l|l|}
\hline
n & decomposition \\
\hline
100 & $C(2)^7\times C(1462491779472195274571694315857495335176880879072)$\\
110 & $C(2)^{11}\times C(8576403641950292891121955131452148838284294200071440)$\\
\hline
\end{tabular}
\end{center}
\end{figure}

%For acknowledgements section, please don't number the section, you need to begin with \section*{Acknowledgements}
\section*{Acknowledgements} The author thanks Andreas Enge for his support on this project, the fruitful discussions we had and a careful reading of this article. We thank Nicolas Th\'{e}riault and all the organizing comitee of the conference CHILE 2009 where the original results of this paper were first presented. We also thank J\'{e}r\^{o}me Milan for his support on issues regarding implementation, especially with the TIFA library.

\medskip

Received June 2009; revised December 2009.

% Please replace the following email addresses by yours.

\medskip
 {\it E-mail address: }biasse@lix.polytechnique.fr\\


\begin{thebibliography}{99}

\bibitem{bach}
    	\newblock E. Bach,
    	\newblock \emph{Explicit bounds for primality testing and related problems},
    	\newblock Mathematics of Computations, \textbf{55} (1990), 335--380.

%\bibitem{siegel}
 %   	\newblock R. Brauer,
  %  	\newblock \emph{On the Zeta-functions of algebraic number fields},
   % 	\newblock American Journal of Mathematics, \textbf{69} (1947), 243--250.

\bibitem{dullmann}
	\newblock J. Buchmann and S. D\"{u}llmann,
	\newblock \emph{On the computation of discrete logarithms in class groups},
	\newblock in ``Advances in Cryptology - CRYPTO '90,"
	\newblock Lecture Notes in Computer Science, \textbf{537} (1991), 134--139.

\bibitem{buchmannProtocol}
    	\newblock J. Buchmann and H.C. Williams,
    	\newblock \emph{A key-exchange system based on imaginary quadratic fields},
    	\newblock Journal of Cryptology, \textbf{1} (1988), 107--118.

\bibitem{Cavallar}
	\newblock S. Cavallar,
	\newblock \emph{Strategies in Filtering in the Number Field Sieve},
	\newblock in ``ANTS-IV: Proceedings of the 4th International Symposium on Algorithmic Number Theory,"
	\newblock Lecture Notes in Computer Science, \textbf{1838} (2000), 209--232.

\bibitem{cohen}
    	\newblock H. Cohen,
    	\newblock ``A course in computational algebraic number theory,"
	\newblock vol 138 of Graduate Texts in Mathematics,
    	\newblock Springer-Verlag, 1991.

\bibitem{contini}
	\newblock S. Contini
	\newblock  ``Factoring integers with the self initializing quadratic sieve,"
	\newblock  Master thesis, University of Georgia, 1997
\bibitem{squfof}
    	\newblock J.E. Gower and S. Wagstaff,
    	\newblock \emph{Square form factorization},
    	\newblock Mathematics of Computations, \textbf{77} (2008), 551--588.

\bibitem{hafner}
    	\newblock J.L. Hafner and K.S. McCurley,
    	\newblock \emph{A rigorous subexponential algorithm for computation of class groups},
    	\newblock J. Amer. Math. Soc., \textbf{2} (1989), 839--850.

\bibitem{JacobsonProtocol}
	\newblock D. H\"{u}hnlein, M.J. Jacobson, S. Paulus and T. Takagi,
	\newblock \emph{A cryptosystem based on non-maximal imaginary quadratic orders with fast decryption},
	\newblock in ``Advances in Cryptology - EUROCRYPT '98,"
	\newblock Lecture Notes in Computer Science, \textbf{1403} (1998), 294--307.

\bibitem{JacobsonPhd}
    	\newblock M. Jacobson,
    	\newblock  ``Subexponential Class Group Computation in Quadratic Orders,"
    	\newblock  Ph.D thesis, Technische Universit\"{a}t Darmstadt, 1999,
 	\newblock Shaker Verlag GmbH.

\bibitem{JacobsonPell}
    	\newblock M.J. Jacobson and H.C. Williams,
% quote the book tile, the first letter of each word should be capitalized.
    	\newblock ``Solving the Pell equation,"
	\newblock  CMS Books in Mathematics,
    	\newblock Springer-Verlag, 2009.

\bibitem{Lenstra2LP}
	\newblock A.K. Lenstra and M.S. Manasse,
	\newblock \emph{Factoring with two large primes (extended abstract)},
	\newblock in ``Advances in Cryptology - EUROCRYPT '90,"
	\newblock Lecture Notes in Computer Science, \textbf{473} (1991), 72--82.

\bibitem{LenstraShanks}
	\newblock A.K. Lenstra,
	\newblock \emph{On the calculation of regulators and class numbers of quadratic fields},
	\newblock in ``Journ\'{e}es arithm\'{e}tiques,"
	\newblock Cambridge Univ. Press, (1982).

\bibitem{Lenstra3LP}
	\newblock P.C. Leyland, A.K. Lenstra, B. Dodson, A. Muffett and S. Wagstaff,
	\newblock \emph{MPQS with Three Large Primes},
	\newblock in ``ANTS-V: Proceedings of the 5th International Symposium on Algorithmic Number Theory,"
	\newblock Lecture Notes in Computer Science, \textbf{2369} (2002), 446--460.

\bibitem{tifa}
	\newblock J. Milan,
	\newblock ``TIFA",
	\newblock http://www.lix.polytechnique.fr/Labo/Jerome.Milan/tifa/tifa.xhtml.

\bibitem{mulders}
	\newblock T. Mulders and A. Storjohann,
	\newblock \emph{Certified linear system solving},
	\newblock Technical report,
	\newblock ETH Z\"{u}rich (2000).

\bibitem{neukirch}
    	\newblock J. Neukirch,
% quote the book tile, the first letter of each word should be capitalized.
    	\newblock ``Algebraic Number Theory,"
	\newblock vol 322 of Comprehensive Studies in Mathematics,
    	\newblock Springer-Verlag, 1999.
	\newblock Translation into english by Norbert Schappacher.

\bibitem{pomerance82}
	\newblock C. Pomerance,
	\newblock \emph{Analysis and comparison of some integer factoring algorithms},
	\newblock in ``Computational methods in number theory I,"
	\newblock Mathematical Centre Tracts, \textbf{154} (1982), 89--139.

\bibitem{Shanks}
	\newblock D. Shanks,
	\newblock \emph{Class number, a theory of factorization, and genera},
	\newblock in ``Proceedings of symposia in pure mathematics,"
	\newblock American Mathematical Society, \textbf{20} (1969), 415--440.

\bibitem{siegel}
	\newblock C.Siegel,
	\newblock \emph{\"{U}ber die Klassenzahl quadratischer Zahlk\"{o}rper},
	\newblock Acta Arithmetica, \textbf{1} (1936), 83--86.

\bibitem{iml}
	\newblock A. Storjohann,
	\newblock ``IML",
	\newblock http://www.cs.uwaterloo.ca/$\sim$astorjoh/iml.html.

\bibitem{Sutherland}
    	\newblock A. Sutherland,
    	\newblock  ``Order Computations in Generic Groups,"
    	\newblock  Ph.D thesis, Massachusetts Institute of Technology, 2007.

\bibitem{vollmer}
	\newblock U. Vollmer,
	\newblock \emph{A note on the {H}ermite basis computation of large integer matrices},
	\newblock in ``ISSAC '03: Proceedings of the 2003 international symposium on Symbolic and algebraic computation,"
	\newblock ACM, (2003), 255--257.

\end{thebibliography}
\end{document}